\documentclass[12pt]{article}
\usepackage{amsthm, amsmath,amssymb}
\usepackage{longtable}
\usepackage[T1]{fontenc}
\usepackage[latin1]{inputenc}
\newtheorem{theorem}{Theorem}[section]
\usepackage{url} 
\usepackage{icomma} 
\usepackage{calc}
\usepackage[hang]{caption}
\usepackage{varioref}
\usepackage{booktabs}
\begin{document}

%
%
%
%

\title{A computer search of maximal partial spreads in $\mathrm{PG}(3, q)$}

\author{Maurizio Iurlo and Sandro Rajola}


\date{\quad}

\maketitle

\begin{abstract}
In this work, by a computer search, new minimum sizes for the maximal partial spreads of $\mathrm{PG}(3,q)$ have been obtained for $q = 8, 9, 16$ and for every $q$ such that $25 \leq q \leq 101$. Furthermore, density results  in the cases $q = 8, 9, 16, 19, 23, 25, 27$ have been obtained. Finally, the known exceptional size $45$ for $q = 7$ has been founded again.
\end{abstract}


\section{Introduction}
A spread of $\mathrm{PG}(3, q)$, projective space of three dimensions  over the
field $\mathrm{GF}(q)$, is a set of mutually skew lines covering the space. A
partial spread is a set of mutually skew lines which  is not a spread. A
partial spread is said to be maximal if it is neither properly contained  in  a
spread or in a partial spread.

Maximal partial spreads have been investigated by several authors, but a
complete knowledge of them is still far.

This work is the natural continuation of the paper ``A new method
to construct maximal partial spreads in $\mathrm{PG}(3, q)$'' \cite{Iura},
where we found new minimums for the sizes of maximal partial spreads of $\mathrm{PG}(3, q)$,
with $q=11, 13, 17, 19, 23$. 

Moreover in \cite{Iura}, for $q=11, 13, 17$, we constructed maximal partial spreads (in the following
Mps) having all the cardinalities between our minimums and those of
the density results found by O. Heden.
In the cases $q=19$ and $q=23$ we did not fill the previous gap,
but we do it here.

In this paper we found new minimums for the sizes of Mps in $\mathrm{PG}(3, q)$,
with $q=8, 9, 16$ and for every $q$ such that  $25\leq q\leq101$. 

Afterwards, we found the necessary cardinalities to fill the gaps between
our minimums and the size $q^{2}-q+2$, and do it for $q=8, 9,  19, 23, 25, 27$. We obtained also 
density results in the case $q=16$.

During the research, we found many known values, such as the exceptional cardinality
$45$ for $q=7$. 

To construct the Mps, we used several programs, written in C language,
and do it by a notebook with processor Intel Core i5-430M, 2.26 GHz,
3 MB L3 cache and 4 GB RAM.

The first program, identified by ``max-intersection program'', it is 
much more efficient version than the one used in \cite{Iura}, and it works
in the following way.

First of all, the program eliminates all the lines meeting some lines
of an initial partial spread. Then it calculates the number of the
remaining lines meeting each remaining line, and  adds to the initial
partial spread the remaining line meeting the maximum number of remaining lines.
The program proceeds in this way until to obtain a Mps.

For all the values of $q$ studied in this paper  we found minimums less than
\[
\left(\left\lceil 2\log_{2}q\right\rceil+1\right)q+1-3q,
\]
where $\left(\left\lceil 2\log_{2}q\right\rceil+1\right)q+1$ is the known minimum for $q$
odd and $q\geq 23$, while for $q$ even the known minimum is much higher
than it.

Furthermore, we used two other versions of the program which again
calculate the number of the lines meeting a fixed line, but select
it when its value is the minimum or the closest to the average. Such programs will
be respectively identified as ``min-intersection program'' and ``middle-intersection program''.

We use these versions to get unknown cardinalities greater than the
found minimums.

Furthermore we write programs which construct several Mps at the
same time. 

For simplicity reasons all the previous programs use the line of plücker
coordinates $\left(0, 0, 0, 0, 0, 1\right)$ as initial line. 

Afterwards, we wrote a program which constructs Mps in the following
way. The program, that we call ``linear program'', 
chooses the first line in the order of construction,
that is the order through which our algorithm constructs the Plücker coordinates 
of the lines,
and eliminates all the lines meeting it. Next, the program chooses
the first of the remaining lines and proceeds similarly until to construct
a Mps. Then the program chooses the second line, in the order of construction,
as first line and constructs the second Mps, and so on. 
So the program constructs $\theta_{3}\theta_{2}/\theta_{1}$
Mps, where $\theta_{r}=q^r+q^{r-1}+\ldots+1$.

The linear program, besides giving many unknown cardinalities, finds
Mps of sizes greater than those obtained by the max-intersection program,
but lower than the previous known minimums.

\section{Our results}
In this work we found new minimums and new density results for the
sizes of Mps of $\mathrm{PG}(3, q)$. In particular, we found new minimums
for $q=8, 9, 16$ and for every $q$ such that $25\leq q\leq 101$. 

Moreover, we found new density results for $q=8, 9, 16, 19, 25, 27$
and the size $149$ for $q=23$, which is the missing value between
the minimum found in \cite{Iura} and the minimum of the density result
found in the same article.
%

Obviously we found many known results, such as the size $45$ for $q=7$,
as already said.

Totally, we constructed about one  million and half Mps or spreads.

The density results found here and which appear in the Table \ref{tab:longtable2} include
also some known values. The knowledge of such values is not specified
for brevity reasons. 
\begin{longtable}{lllll}
\caption{New sizes of maximal partial spreads in $\mathrm{PG}(3, q)$}
\label{tab:longtable2}\\
\toprule
$q$ & Min. & Previous min. & Density results & Previous density results \\ 
\hline 
\endhead
\endfoot
8 & 30 & 41 & 31\thinspace--\thinspace55 & 56\thinspace--\thinspace58 \tabularnewline \hline 
9 & 36 & 46 & 37\thinspace--\thinspace45 & 46\thinspace--\thinspace74\tabularnewline \hline
16 & 87 & 145 & 88\thinspace--\thinspace221, 225\thinspace--\thinspace231  & 240\thinspace--\thinspace242 \tabularnewline \hline
19 &  & 114 & 147\thinspace--\thinspace181 & 115\thinspace--\thinspace146; 182\thinspace--\thinspace344\tabularnewline \hline
23 &  & 148 & 149 & 150\thinspace--\thinspace508;\tabularnewline \hline
25 & 173 & 276 & 174\thinspace--\thinspace313 & 314\thinspace--\thinspace602\tabularnewline \hline
27 & 193 & 298 & 194\thinspace--\thinspace367 & 368\thinspace--\thinspace704\tabularnewline \hline
29 & 210 & 320 &  & \tabularnewline \hline
31 & 231 & 342 &  & \tabularnewline \hline
32 & 238 & 545 &  & \tabularnewline \hline
37 & 306 & 445 &  & \tabularnewline \hline
41 & 345 & 493 &  & \tabularnewline \hline
43 & 372 & 517 &  & \tabularnewline \hline
47 & 417 & 612 &  & \tabularnewline \hline
49 & 474 & 638 &  & \tabularnewline \hline
53 & 488 & 690 &  & \tabularnewline \hline
59 & 569 & 768 &  & \tabularnewline \hline
61 & 600 & 794 &  & \tabularnewline \hline
64 & 623 & 1665 &  & \tabularnewline \hline
67 & 672 & 939 &  & \tabularnewline \hline
71 & 732 & 995 &  & \tabularnewline \hline
73 & 761 & 1023 &  & \tabularnewline \hline
79 & 848 & 1107 &  & \tabularnewline \hline
81 & 873 & 1135 &  & \tabularnewline \hline
83 & 903 & 1163 &  & \tabularnewline \hline
89 & 968 & 1247 &  & \tabularnewline \hline
97 & 1102 & 1456 &  & \tabularnewline \hline
101 & 1160 & 1516 &  & \\
\bottomrule
\end{longtable}
From the already known results and from our results, we get the following theorem.

\begin{theorem}
In $\mathrm{PG}(3, q)$, for every $q$ such that $5 \leq q \leq 101$, there are maximal partial spreads of size less than
\[
\left(\left\lceil 2\log_{2}q\right\rceil+1\right)q+1-3q.
\]
\end{theorem}
Concerning density results, in the case $q=16$ we haven't found all the unknown cardinalities included between the minimum value we found and the biggest unknown cardinality, in spite of numerous attempts. This is really unexpected, because in the other cases we have found all the unknown cardinalities in a very easy way.

In addition, from the already known results and from our results, we get the following theorem.
\begin{theorem}
In $\mathrm{PG}(3, q)$, for every $q$ such that $5 \leq q \leq 27$ and $q\neq 16$, there is a maximal partial spread of size $n$
for any integer $n$ in the interval
\[
\left(\left\lceil 2\log_{2}q\right\rceil+1\right)q+1-3q \leq n \leq q^2-q+2.
\]
\end{theorem}
For every new example of Mps, we specify the program through
which we have obtained it. 
Obviously, we  have obtained several results using
 different programs. 
 

We give some examples about the execution time of the programs.

For $q=7$ the linear program finds all the sizes between $27$
and $45$, and does it in $1,37$ seconds.

For $q=8$ the linear program constructs, in $5,95$ seconds,
$4096$ Mps or spreads having all the cardinalities between $33$ and
$52$, and the cardinalities $54$, $56$, $57$ and $65$.

For $q=9$ the linear program constructs, in $16,89$ seconds,
$7462$ Mps or spreads having all the cardinalities between $41$ and
$69$, and the cardinalities $71$, $72$ and $82$.


The max-intersection program  gives, for $q=8$, the cardinality $30$ in $0,46$ seconds; for $q=9$ the cardinality $36$ in $0,87$ seconds; for $q=16$ the cardinality $87$ in $19,80$ seconds and, for $q=32$, the cardinality $238$ in $648,09$ seconds.


For $q=71$ the max-intersection program  gives the cardinality $732$ in $2571,34$ seconds, the cardinality $785$ in $119,78$ seconds and the cardinality $983$ (that is lower than the previous known minimum) in $42,71$ seconds.

\section{Some new examples of maximal partial spreads}
In this number we report the plücker coordinates of the lines of some Mps
that we find.\newline
For every reported Mps, we firstly write the plücker coordinates of the lines of the initial partial spread, and then the order numbers $i$ of the added lines, whose plücker coordinates can be determined through the formulas: 
\begin{align*}
p_{01}&=1,\\
p_{02}&=i\operatorname{mod}q,\\
p_{03}&=\left\lfloor i/q\right\rfloor \operatorname{mod}q,\\
p_{12}&=\left\lfloor i/q^{2}\right\rfloor \operatorname{mod}q,\\
p_{13}&=\left\lfloor i/q^{3}\right\rfloor \operatorname{mod}q,\\
p_{23}&=(p_{02}p_{13}-p_{03}p_{12}) \operatorname{mod}q.
\end{align*}
\quad \newline
Maximal partial spread of size $30$ for $q=8$. \newline
\quad\newline
Initial lines:\newline
(0, 0, 0, 0, 0, 1), (1, 4, 1, 0, 6, 5), (1, 0, 0, 1, 6, 0), (1, 1, 2, 2, 6, 2), (1, 1, 3, 3, 6, 3), (1, 1, 4, 4, 6, 0), (1, 1, 5, 5, 6, 1), (1, 1, 6, 6, 6, 4).\newline
Added lines:\newline
24, 2367, 231, 3708, 455, 2394, 3784, 1165, 180, 3971, 2134, 2589, 1893, 1808, 631, 3883, 382, 1462, 2063, 708, 810, 1537.\newline
\quad \newline
\noindent Maximal partial spread of size $210$ for $q=29$. In order to construct this Mps, we choose sixtyone lines from the Bruen-Hirschfeld's spread.\newline
\quad \newline
\noindent Initial lines:\newline
(0, 0, 0, 0, 0, 1), (9, 0, 9, -1, 0, 1), (16, 8, 12, 0, 2, 1), (20, 28, 21, 3, 4, 1), (28, 14, 7, 8, 6, 1), (13, 7, 28, 15, 8, 1), (28, 0, 7, -4, 0, 1), (24, 6, 19, 0, 4, 1), (1, 21, 26, 12, 8, 1), (13, 25, 28, 3, 12, 1), (5, 27, 25, 2, 16, 1), (4, 0, 23, -9, 0, 1), (20, 13, 21, 0, 6, 1), (25, 2, 15, 27, 12, 1), (6, 1, 5, 14, 18, 1), (9, 15, 20, 19, 24, 1), (9, 0, -3, 3, 0, 1),(16, 8, 0, 4, 2, 1), (20, 28, 9, 7, 4, 1), (28, 14, 24, 12, 6, 1), (13, 7, 16, 19, 8, 1), (28, 0, -12, 12, 0, 1), (24, 6, 0, 16, 4, 1), (1, 21, 7, 28, 8, 1), (13, 25, 9, 19, 12, 1), (5, 27, 6, 18, 16, 1), (4, 0, -27, 27, 0, 1), (20, 13, 0, 7, 6, 1), (25, 2, 23, 5, 12, 1), (6, 1, 13, 21, 18, 1), (9, 15, 28, 26, 24, 1), (1, 2, 3, 1, 2, 1), (16, 16, 12, 4, 4, 1), (23, 25, 27, 9, 6, 1), (24, 12, 19, 16, 8, 1), (16, 18, 17, 25, 10, 1), (20, 26, 21, 7, 12, 1), (23, 19, 2, 20, 14, 1), (7, 9, 18, 6, 16, 1), (7, 8, 11, 23, 18, 1), (24, 28, 10, 13, 20, 1), (25, 23, 15, 5, 22, 1), (1, 5, 26, 28, 24, 1), (25, 15, 14, 24, 26, 1), (20, 7, 8, 22, 28, 1), (20, 22, 8, 22, 1, 1), (25, 14, 14, 24, 3, 1), (1, 24, 26, 28, 5, 1), (25, 6, 15, 5, 7, 1), (24, 1, 10, 13, 9, 1), (7, 21, 11, 23, 11, 1), (7, 20, 18, 6, 13, 1), (23, 10, 2, 20, 15, 1), (20, 3, 21, 7, 17, 1), (16, 11, 17, 25, 19, 1), (24, 17, 19, 16, 21, 1), (23, 4, 27, 9, 23, 1), (16, 13, 12, 4, 25, 1), (1, 27, 3, 1, 27, 1).\newline
\quad\newline
Added lines:\newline
677253, 504585, 521560, 449301, 625597, 347945, 489072, 36563, 240119, 509323, 616226, 330155, 82544, 121871, 174971, 187236, 138497, 157222, 346096, 108275, 124884, 147268, 601567, 391027, 429148, 109152, 145311, 432435, 550591, 697973, 25022, 191857, 173609, 589158, 459617, 129059, 206486, 596160, 367651, 56530, 337034, 658419, 317597, 55325, 603163, 52495, 107491, 451648, 683065, 148086, 285155, 116416, 302602, 337486, 150168, 477206, 196604, 506753, 274083, 561501, 33049, 42382, 458736, 70067, 569409, 441523, 416479, 80220, 243346, 537506, 516647, 89547, 328090, 212003, 98520, 109483, 234264, 215347, 245551, 503476, 528854, 21953, 385516, 271778, 527360, 189087, 423060, 232916, 38771, 286659, 112330, 669444, 296968, 277363, 182475, 583897, 482186, 160767, 110259, 38321, 642746, 341987, 105983, 122833, 49273, 531042, 304797, 519957, 115948, 644653, 594328, 395375, 650790, 492067, 662581, 113012, 494299, 7416, 498804, 103763, 167220, 272780, 58186, 47265, 200268, 372443, 421720, 605597, 76597, 464438, 706631, 90926, 437079, 453946, 510541, 27980, 70865, 152474, 344471, 410036, 27349, 111830, 197156, 197918, 202937, 241613, 254354, 370916, 379262, 397943.
677253, 504585, 521560, 449301, 625597, 347945, 489072, 36563, 240119, 509323, 616226, 330155, 82544, 121871, 174971, 187236, 138497, 157222, 346096, 108275, 124884, 147268, 601567, 391027, 429148, 109152, 145311, 432435, 550591, 697973, 25022, 191857, 173609, 589158, 459617, 129059, 206486, 596160, 367651, 56530, 337034, 658419, 317597, 55325, 603163, 52495, 107491, 451648, 683065, 148086, 285155, 116416, 302602, 337486, 150168, 477206, 196604, 506753, 274083, 561501, 33049, 42382, 458736, 70067, 569409, 441523, 416479, 80220, 243346, 537506, 516647, 89547, 328090, 212003, 98520, 109483, 234264, 215347, 245551, 503476, 528854, 21953, 385516, 271778, 527360, 189087, 423060, 232916, 38771, 286659, 112330, 669444, 296968, 277363, 182475, 583897, 482186, 160767, 110259, 38321, 642746, 341987, 105983, 122833, 49273, 531042, 304797, 519957, 115948, 644653, 594328, 395375, 650790, 492067, 662581, 113012, 494299, 7416, 498804, 103763, 167220, 272780, 58186, 47265, 200268, 372443, 421720, 605597, 76597, 464438, 706631, 90926, 437079, 453946, 510541, 27980, 70865, 152474, 344471, 410036, 27349, 111830, 197156, 197918, 202937, 241613, 254354, 370916, 379262, 397943.

\section{Conclusion}
This work has the aim not only of finding new minimum sizes for the
maximal partial spreads of $\mathrm{PG}(3,q)$, but also of giving,
as an obvious consequence, a theoretical indication and therefore
a new impulse to the research, which  stopped  seven years ago.
In fact, the last results go back to the year 2003, when A.\thinspace G\'{a}cs and T.\thinspace Sz\H{o}nyi managed to lower the previous minimums remarkably.\\
However, the gaps between the Glynn's lower bound and the minimums
we know until now still appeared much too large.\\
Here, for the values of $q$ that we study, we succeed in getting
a reduction up to  $70\%$ of the previous gaps, as happens
in the case $q=64$. \\
Moreover, we have noted not only that the new minimums are quite
lower than the previous one, but also that an essential difference
between the cases $q$ even and $q$ odd does not appear. Only the case $q=16$ has been different from the others, but only for the density results.\\
However, it is possible to develop the computer search, too. We are putting right a new program ourselves, through which you
can investigate values of $q$ much larger than those studied in this
paper.

\quad\\
\begin{center}
--------------------
\end{center}
\quad\\

\noindent Maurizio Iurlo\\
\noindent Largo dell'Olgiata, 15 \\
\noindent 00123 Roma \\
\noindent Italy\\
\noindent maurizio.iurlo@istruzione.it\\
\noindent www.maurizioiurlo.com \\
\quad\\
\noindent Sandro Rajola\\
\noindent Istituto Tecnico per il Turismo ``C. Colombo''\\
\noindent Via Panisperna, 255\\
\noindent 00184 Roma\\
\noindent Italy\\
\noindent sandro.rajola@istruzione.it


\begin{thebibliography}{999}

\bibitem {BA} J.\thinspace B\'{a}rat, A.\thinspace Del Fra, S.\thinspace Innamorati and L.\thinspace Storme, Minimal Blocking Sets in $\mathrm{PG}(2,8)$ and Maximal Partial Spreads in $\mathrm{PG}(3,8)$, 
\emph{Des. Codes Cryptogr.} \textbf{31} (2004), \mbox{15--26}.

\bibitem {BEU} A.\thinspace Beutelspacher, Blocking sets and partial spreads in finite projective spaces,
\emph{Geom. Dedicata} \textbf{9} (1980), \mbox{425--449}.

\bibitem {BL} A.\thinspace Blokhuis, Note on the size of a blocking set in $\mathrm{PG}(2,p)$, \emph{Combinatorica} \textbf{14} (1994), \mbox{111--114}. 

\bibitem {B} A.\thinspace A.\thinspace Bruen, Partial spreads and replaceable
nets, \emph{Canad. J. Math.} \textbf{23 }(1971), \mbox{381--392}.

\bibitem {BH} A.\thinspace A.\thinspace Bruen and J.\thinspace W.\thinspace P.\thinspace
Hirschfeld, Applications of line geometry over finite fields. I:
The twisted cubic, \emph{Geom. Dedicata} \textbf{6 }(1977), \mbox{495--509}.

\bibitem {BT} A.\thinspace A.\thinspace Bruen and J.\thinspace A.\thinspace Thas,
Partial spreads, packings and hermitian manifolds in $\mathrm{PG}(3,q)$, \emph{Math.
Z.} \textbf{151 }(1976), \mbox{207--214}.

\bibitem {F} J.\thinspace W.\thinspace Freeman, Reguli and pseudo-reguli in $\mathrm{PG}(3,s^2)$, \emph{Geom. Dedicata} \textbf{9} (1980), \mbox{267--280}.

\bibitem {Gal} F.\thinspace Fuhlendorf and A.\thinspace Stuht,  Software project Calculator for finite Galois fields, http://www.fh-wedel.de/mitarbeiter/iw/eng/r-d/done/sw-projects/galoisfield/

\bibitem {G} D.\thinspace G.\thinspace Glynn, A lower bound for maximal
partial spreads in  $\mathrm{PG}(3,q)$, \emph{Ars Comb.} \textbf{13}
(1982), \mbox{39--40}.

\bibitem {GS} A.\thinspace G\'{a}cs and T.\thinspace Sz\H{o}nyi, On maximal partial spreads
in $\mathrm{PG}(n,q)$, \emph{Des. Codes Cryptogr.} \textbf{29 } (2003), \mbox{123--129}.

\bibitem {H1} O.\thinspace Heden,   A greedy search for maximal partial spreads
in $\mathrm{PG}(3,7)$, \emph{Ars Comb.} \textbf{32} (1991), \mbox{253--255}.

\bibitem {H2} O.\thinspace Heden, Maximal partial spreads and the modular
$n$-queen problem, \emph{Discrete Math.} \textbf{120 }(1993), \mbox{75--91}.

\bibitem {H3}  O.\thinspace Heden, Maximal partial spreads and the modular
$n$-queen problem. II, \emph{Discrete Math.} \textbf{142 }(1995), \mbox{97--106}.

\bibitem {H4}  O.\thinspace Heden, Maximal partial spreads in $\mathrm{PG}(3,5)$, \emph{Ars
Comb.} \textbf{57 }(2000), \mbox{97--101}.

\bibitem {H45} O.\thinspace Heden, A Maximal Partial Spread of Size $45$ in $\mathrm{PG}(3,7)$, \emph{Des. Codes Cryptogr.} \textbf{22} (2001), \mbox{331--334}. 

\bibitem {H5}  O.\thinspace Heden, Maximal partial spreads and the modular
$n$-queen problem. III, \emph{Discrete Math.} \textbf{243 }(2002), \mbox{135--150}.

\bibitem {H6} O.\thinspace Heden, No maximal partial spread of size $115$ in $\mathrm{PG}(3,11)$, \emph{Ars Comb}. \textbf{66} (2003).

\bibitem {HAL1} O.\thinspace Heden, S.\thinspace Marcugini, F.\thinspace Pambianco and L.\thinspace Storme,  On the non-existence of a maximal partial spread of size $76$ in $\mathrm{PG}(3,9)$, \emph{Ars Comb.} \textbf{89} (2008).

\bibitem {Hir}  J.\thinspace W.\thinspace P.\thinspace Hirschfeld, \emph{Finite projective spaces of three dimension}, Oxford University Press, Oxford, 1985.

\bibitem {Iura}  M.\thinspace Iurlo and S.\thinspace Rajola, \emph{A new method to construct maximal partial spreads of smallest sizes in} $\mathrm{PG}(3, q)$. In: Error-Correcting Codes, Cryptography and Finite Geometries, Contemporary Mathematics, vol. 523, Amer. Math. Soc., Providence, RI, 2010, \mbox{89--107}.

\bibitem {JS} D.\thinspace Jungnickel and L.\thinspace Storme, A note on maximal partial spreads with deficiency $q+1$, $q$ even, \emph{J. Combin. Theory Ser. A} \textbf{102} (2003), \mbox{443--446}.

\bibitem {M} D.\thinspace M.\thinspace Mesner, Sets of disjoint lines in $\mathrm{PG}(3,q)$, \emph{Canad. J. Math.} \textbf{19 }(1967), \mbox{273--280}.

\bibitem {RAJ}  S.\thinspace Rajola and M.\thinspace S.\thinspace Tallini, 
Maximal partial spreads in  $\mathrm{PG}(3,q)$, \emph{ J. Geom} \textbf{85} (2006), \mbox{138--148}.

\bibitem {R}  S.\thinspace Rajola, A construction, in $\mathrm{AG}(2,q)$, of
maximal partial spreads of  $\mathrm{PG}(3,q)$, submitted to \emph{New Zealand J. Math.}

\bibitem {SOI} L.\thinspace H.\thinspace Soicher, Computation of partial spreads, \newline http://www.maths.qmul.ac.uk/$\sim$leonard/partialspreads/

\end{thebibliography}
\end{document}